\documentstyle[draft,amscd,syntonly,amssymb]{amsart}

\theoremstyle{plain}
\newtheorem{Thm}{Theorem}
\newtheorem{Prop}[Thm]{Proposition}
\newtheorem{Cor}[Thm]{Corollary}
\newtheorem{Lem}[Thm]{Lemma}

 \theoremstyle{definition}
\theoremstyle{remark}

\numberwithin{equation}{section}

\begin{document} 
 \title{ HAMILTONIAN SYMPLECTOMORPHISMS AND THE BERRY PHASE} 
    
 \author{ ANDR\'{E}S   VI\~{N}A}
\address{Departamento de F\'{i}sica. Universidad de Oviedo.   Avda Calvo
 Sotelo.     33007 Oviedo. Spain. } 
 \email{cuevas@@pinon.ccu.uniovi.es}
\thanks{The author was partially supported by DGES, grant DGES-96PB0538}
  \keywords{Group of symplectomorphisms, Geometric quantization, Berry phase, Coadjoint orbits}
 \subjclass{Primary: 53C15; Secondary: 53C12, 81S10}

 \maketitle
\begin{abstract}

On the space ${\cal L}$, of loops in    the group of Hamiltonian symplectomorphisms 
of a  symplectic quantizable  manifold,  we define 
a closed ${\bf Z}$-valued  $1$-form $\Omega$. If $\Omega$ vanishes,
  the prequantization map can be  extended to a group representation. On ${\cal L}$ 
one can  define 
an action integral as an ${\bf R}/{\bf Z}$-valued function, and 
the cohomology class $[\Omega]$  is the obstruction to the   lifting
of that action integral to an ${\bf R}$-valued function. The form $\Omega$ also defines
a natural grading on $\pi_1({\cal L})$.

 \end{abstract}

 \section {Introduction} \label{S:intro}
In the process of quantization of a symplectic manifold $(M,\omega)$ it is
 necessary to fix a polarization $I$, then the corresponding quantization ${\cal Q}_I$
is the space  of the sections of a prequantum bundle $L$, which are parallel
along the leaves of the polarization $I$ \cite{nW92}. 
The identification of the ${\cal Q}_I$ obtained by fixing different polarizations is one of the goals of the
geometric quantization, but ``the theory is far from achieving this goal"
 \cite[page 267]{vG84}. 
This issue has been treated in several particular cases: 
The identification of the quantizations of the moduli space
of flat connections on a closed surface has been studied in \cite{sA91} and in
\cite{nH90}; the case when $M$ is a torus has been treated in \cite{aV96}. The problems 
 involved in an
identification of the spaces ${\cal Q}_I$  
 were analysed in \cite{aV00}, when the polarizations considered are of type K\"ahler.

Here we consider a similar situation: If $\{\psi_t\,|\,t\in[0,1]\}$ is a Hamiltonian isotopy of $M$ 
 \cite{M-S98} and $F$ a foliation on $M$,
the  action of $\psi_t$ produces a family $F_t$ of foliations. 
We have the spaces ${\cal Q}_{F_t}$ of sections of $L$ which are ``polarized" with respect to $F_t$; i.e.,
sections parallel along the leaves of $F_t$. 
We shall construct
 isomorphisms $\tau\in{\cal Q}_F\rightarrow \tau_t\in{\cal Q}_{F_t}$, which permit us
 to ``transport"
the vectors in ${\cal Q}_F$ to the spaces ${\cal Q}_{F_t}$ in a continuous way.
In general this transport has non vanishing ``curvature"; that is, it depends
 on the isotopy which joins
a given symplectomorphism with $\text{id}$.

In the prequantization process of $(M,\omega)$ one assigns to each function $f$ on $M$ an operator
 ${\cal P}_f$  \cite[p. 57-59]{jS80}, which acts on the space $\Gamma(L)$ of sections of $L$.
The map ${\cal P}$ is a representation of $\text{Vect}_H(M)$, the algebra of Hamiltonian vector fields
on $M$. There are obstructions to extend this representation to a representation of
$\text{Ham}(M)$, the
group of Hamiltonian symplectomorphisms of $M$ \cite{M-S98}. 
We analyse the relation between
these obstructions and the curvature of the aforementioned transport.

 If the Hamiltonian isotopy  $\psi_t$ is a loop
in the group $\text{Ham}(M)$ and $N$ is a Lagrangian 
leaf of the foliation $F$, then $\psi_t(N)$
is a loop  of submanifolds of $M$ and the corresponding Berry phase is defined
\cite{aW90}. We prove the existence of
a  number $\kappa(\psi)\in U(1)$, which depends only on the loop $\psi$, 
and that relates {\em any section} $\rho$ with $\rho_1$, the section 
resulting of the transport of $\rho$,
by the formula
 $\rho_1=\kappa(\psi)\rho$.  So $\kappa(\psi)$ is the ``holonomy" of the transport along $\psi$.
 It turns out that the holonomy of our transport is essentially the
Berry phase of the loop $\psi_t(N)$.   Using 
the map $\kappa$ we construct on ${\cal L}$, 
the space of loops in $\text{Ham}(M)$ based at $\text{id}$,
 a closed $1$-form $\Omega$.
 The vanishing of $\Omega$ is equivalent to the invariance of the Berry 
phase under deformations of the loop $\psi$. We will prove that 
there is a well-defined an ${\Bbb R}/{\Bbb Z}$-valued action integral on ${\cal L}$.
The exactness of $\Omega$ is  equivalent to the 
existence of a lift of the  action integral to an ${\Bbb R}$-valued
map.   
The integral   of the form $\Omega$ along   a loop $\phi^s$ in ${\cal L}$ is 
in fact the winding number of
the map $s\in S^1\mapsto\kappa(\phi^s)\in U(1)$, so $\Omega$ is ${\Bbb Z}$-valued. 
This property permits to define a grading on
$\pi_2(\text{Ham}(M))$ compatible with the group structure.

In Section 2 is introduced the transport of vectors $\tau\in {\cal Q}_F$ to vectors
  $\tau_t\in{\cal Q}_{F_t}$. Such a transport is determined by the differential equation which
it generates; that is, 
$$\frac{d\tau_t}{dt}=\zeta(F_t,\tau_t),$$ 
where $\zeta$ is a section of $L$. The condition $\tau_t\in {\cal Q}_{F_t}$
gives rise to an equation for  $\zeta$. This equation does not determine uniquely 
$\zeta$; however 
it is possible to  choose a {\em natural} solution for $\zeta$
using the time dependent Hamiltonian $f_t$ which generates the isotopy.
If the isotopy is closed, i.e. $\psi_1=\text{id}$, given a leaf $N$ of $F$, it
 is easy to show the
existence of a constant $\kappa$ such that $\tau_{1|N}=\kappa\tau_{|N}$ for all 
$\tau\in{\cal Q}_F$. So one can define the holonomy for the transport of such
sections  $\tau_{|N}$.
 In this Section we also study the relation between this
holonomy and the Berry phase of the loop $\psi_t(N)$ of Lagrangian submanifolds
of $M$.

The Section 3 is concerned with the properties of $\kappa(\psi)$. 
First we prove its existence and determine its expression in terms of the 
Hamiltonian function and the symplectic form. 
Given 
$\{\psi_t\,|\, t\in[0,1]\}$ a loop in $\text{Ham}(M)$, if $q$
is a point of $M$, then the {\em general} action integral around the closed 
curve $\psi_t(q)$ is $\int _S\omega$,
where $S$ is any $2$-submanifold bounded by the curve $\psi_t(q)$. However, for these
particular curves one can also define the action integral
\begin{equation}\label{exp(2pii}
{\cal A}(\psi(q))=\int_S\omega-\int_0^1f_{t}(\psi_t(q))\,dt.
\end{equation}
${\cal A}(\psi(q))$ is well-defined considered as  an element of ${\Bbb R}/{\Bbb Z}$.
Using   results of Section 2 about the transport of polarized sections we will
prove that ${\cal A}(\psi(q))$ 
 is {\em independent} of the point $q\in M$ and that 
$\kappa(\psi)=\text{exp}(2\pi i{\cal A}(\psi)).$ 

 In  \cite{aW89} Weinstein   defined a representation ${\bf A}$ of $\pi_1(\text{Sym}(M))$
as follows: ${\bf A}(\psi)$ is 
the mean value over $q$ of the general action integrals around the curves
$\psi_t(q)$. When $\psi_t$ is a $1$-parameter subgroup
generated by a Hamiltonian function $f$, ${\bf A}$ and $\kappa$    are related by
$\kappa(\psi)=\text{exp}(2\pi i{\bf A}(\psi))$, assumed that the Hamiltonian
function $f$ is normalized
so that $\int f\omega^n=0$. The domain of the map $\kappa$ is less general than the
domain of ${\bf A}$; however the restriction 
to the Hamiltonian symplectomorphisms allows us
 to introduce the second   summand  in (\ref{exp(2pii}), so we obtain
an invariant without averaging on $M$; that is, in contrast with ${\bf A}(\psi)$ the  
value $\kappa(\psi)$ can be calculated pointwise.    
$\kappa$ is not invariant under homotopies; this fact has an interesting meaning.
One can define a $1$-form on ${\cal L}$ as follows:
Given a curve $\psi^s$ in ${\cal L}$, and denoting by $Z$ the vector field defined 
by this curve, the action of the $1$-form $\Omega$ on $Z$ is given by
$$\Omega(Z)=\frac{-1}{2\pi i}\frac{d}{ds}\big( \log(\kappa(\psi^s)) \big).$$
Hence the vanishing of $\Omega$ is equivalent to the invariance of $\kappa(\psi)$ with respect
to deformations of the isotopy $\psi$. The property   $\Omega=0$ is also a
 sufficient condition for ${\cal P}$ extends to a representation of
$\tilde{\text{Ham}}(M)$, the universal cover of $\text{Ham}(M)$.

In Section 4 we prove that $\Omega$ is a closed $1$-form that defines an element of
$H^1({\cal L},\,{\Bbb Z})$. We will also find a simple interpretation
of the cohomology class of $\Omega$; it is the obstruction for the lifting of ${\cal A}$ 
  to an ${\Bbb R}$-valued map. The identification
of $\pi_1({\cal L})$ with $\pi_2(\text{Ham}(M))$ will allow us to
define a grading on the group $\pi_2(\text{Ham}(M))$ by means of the form
$\Omega$.

In Section 5  we consider as symplectic manifold a coadjoint orbit 
 of the group $SU(2)$. There are orbits ${\cal O}$
  diffeomorphic to $S^2$ and
for these manifolds it is easy to determine the value of $\kappa$ on
the loops which are $1$-parameter subgroups in
$\text{Ham}({\cal O})$. With this example we  check
the general properties of $\kappa$ stated in Section 3.

  \section{Loops of  submanifolds and the Berry phase}

Let \( M \) be a connected, compact, symplectic \(C^{\infty} \) manifold of dimension
\(2n \), with symplectic form   \( \omega \). 
 Let us suppose that \((M,\omega)\) is {\em quantizable}, in other words,
 we assume that  
\(\omega\) defines a cohomology class
 in \(H^2(M, \Bbb R)\) which belongs 
to the image of \(H^2(M, \Bbb Z)\) in \(H^2(M, \Bbb R)\) \cite [page  158]{nW92}. 
Then there exists a smooth 
Hermitian line bundle on \(M\) whose first Chern class is \([\omega]\), and on 
this bundle is defined a
connection $D$ compatible with the Hermitian structure and whose curvature is
 \(-2\pi i\omega\). The bundle and the connection are not 
uniquely determined by \(\omega\). 
The family of all possible pairs (line bundle, connection) can be labelled by the elements
of \(H^1(M, U(1))\)  \cite [page  161]{nW92}. From now on we suppose
 that a ``prequantum bundle" \(L\) and a connection 
\(D\) have been fixed, unless it is otherwise indicated.

\smallskip

Let $F$ be a foliation on $M$. If $\tau$ is a 
   $C^{\infty}$  section of $L$    such that $D_A\tau=0$, for all $A\in F$, then $\tau$
is called an $F$-{\em polarized} section, and  the space of $F$-polarized sections of $L$
is denoted by ${\cal Q}_F$. 
 
  Let $\{\psi_t\,|\, t\in[0,1]\}$ be the Hamiltonian isotopy in $M$ generated by the time
 dependent Hamiltonian function $f_t$.
That is,
$$\frac{d\psi_t}{dt}=X_t\circ\psi_t,\,\,\,\iota_{X_t}\omega=-df_t,\,\,\, \psi_0=\text{id}.$$
Then for each $t$ we have a distribution ${F}_t:=(\psi_t)_*(F)$. Moreover, if $N$
is an integral submanifold of $F$ then $N_t:=\psi_t(N)$ is an  
integral submanifold of ${F}_t$. The family $N_t$ is  an isodrastic 
deformation of $N$ \cite{aW90}.

\smallskip

Given $\tau$ an $F$-polarized section of $L$, we
want to define a continuous family $\tau_t$ of sections of $L$ such
 that $\tau_0=\tau$ and
$\tau_t$ is $F_t$-polarized, 
for all $t$. The continuity condition means that there is a section $\zeta$
of $L$ such that
\begin{equation}\label{continuity}
\tau_{t+s}=\tau_t+s\zeta(\tau_t) +O(s^2),
\end{equation}
where $O(s^2)$ is relative to the uniform $C^1$-norm   in the space  $\Gamma(M,L)$
of $C^{\infty}$ sections of $L$.   We will see the 
restrictions on $\zeta$ involved by the
continuity condition (\ref{continuity}), but first of all
we start with a previous result.

Given the isotopy $\psi_t$, each section $\rho$ of $L$ determines
a family $\rho_t$ of sections  by the equation 
\begin{equation}\label{transrho}
\frac{d\rho_t}{dt}=-D_{X_t}\rho_t-2\pi if_t\rho_t,\,\,\,\,\, \rho_0=\rho.
\end{equation}

\begin{Prop}\label{ProposiNueva}
Let $A$ be a vector field on $M$.
If the family $\rho_t$ of sections of $L$ satisfies 
(\ref{transrho}), then
$D_A\rho=0$ implies
  $D_{A_t}\rho_t=0$, for $A_t=(\psi_t)_*(A)$.
\end{Prop}

\begin{pf}
For a fixed $t$ one has 
$$\Big(\frac{d\psi_{t'}(q)}{dt'}\Big)_{|t'=t}=X_t(\psi_t(q)).$$
If we put $t'=t+s$ and $\phi_s:=\psi_{t+s}\circ\psi_t^{-1}$, then
\begin{equation}\label{phis}
\Big(\frac{d\phi_s(p)}{ds}\Big)_{|s=0}=X_t(p).
\end{equation}
 As $\{\phi_s\}$ satisfies (\ref{phis}), 
then for  
 the vector field $A'_s=(\phi_s)_*(A_t)$ we have
\begin{equation}\label{corche}
A'_s=A_t-s[X_t,A_t]+O(s^2).
\end{equation}
Since $A_{t+s}=A'_s$ one has
\begin{equation}\label{Apuntot}
\Dot A_t:=\Big(\frac{dA_{t'}}{dt'}\Big)_{|t'=t}=\Big(\frac{d A'_s}{ds} \Big)_{|s=0}=-[X_t,A_t].
\end{equation}

On the other hand
 \begin{equation}\label{d/dt(AAtrho}
\frac{d}{dt}\big( D_{A_t}\rho_t \big)=D_{\Dot A_t}\rho_t+
D_{A_t}\big( -D_{X_t}\rho_t-2\pi if_t\rho_t   \big).
\end{equation}
As the curvature of $D$ is $-2\pi i\omega$
\begin{equation}\label{auxi-DAt}
-D_{A_t}D_{X_t}\rho_t-D_{[X_t,A_t]}\rho_t=-2\pi i\omega(X_t,A_t)\rho_t-D_{X_t}D_{A_t}\rho_t.
\end{equation}
Since $\iota_{X_t}\omega=-df_t$, by (\ref{Apuntot}) from (\ref{d/dt(AAtrho}) and
(\ref{auxi-DAt}) it follows
\begin{equation}
\frac{d}{dt}\big(D_{A_t}\rho_t  \big)=-D_{X_t}(D_{A_t}\rho_t)-2\pi if_tD_{A_t}\rho_t.
\end{equation}
This is a first order differential equation for the section $\xi(t):=D_{A_t}\rho_t$; if $\xi(0)=D_A\rho$ 
is zero, then $D_{A_t}\rho_t=0$ for every $t$  by the uniqueness of  solutions.
\end{pf}

Given $\rho\in\Gamma(L)$, the family $\{\rho_t\}$ which satisfies the equation
(\ref{transrho}) defines a ``transport" of $\rho$, ``along the isotopy" $\psi=\{\psi_t\}$,
with the property that $\rho_t\in {\cal Q}_{F_t}$ if $\rho\in {\cal Q}_F$. 
The time-$1$ section $\rho_1$ will be 
denoted by ${\cal T}_{\psi}(\rho)$.

 By   Proposition \ref{ProposiNueva}  as section $\zeta$ in (\ref{continuity}) can be taken 
\begin{equation}\label{zetataut}
\zeta(\tau_t)=-D_{X_t}\tau_t -2\pi if_t\tau_t.
\end{equation} 
In general this is not the unique possibility for $\zeta$. In fact if $A$ is a vector field with
$A_p\in F(p)\subset T_pM$, using (\ref{corche}) and (\ref{continuity}) one has
$$D_{A'_s}\tau_{t+s}=D_{A_t}\tau_t+s\big(D_{A_t}\zeta-D_{[X_t,A_t]}\tau_t \big)+O(s^2).$$
As $A_s'=(\psi_{t+s})_*(A)\in F_{t+s}$, the conditions $D_{A'_s}\tau_{s+t}=0$, 
and $D_{A_t}\tau_t=0$ imply
\begin{equation}\label{transport}
D_{A_t}\zeta=D_{[X_t,A_t]}\tau_t\,\,\,\,\text{for every}\,\,\, A_t\in F_t.
\end{equation}
This is the equation for $\zeta$; and it is straightforward 
to check that the $\zeta$ defined in (\ref{zetataut}) satisfies
(\ref{transport}).

 The   solution (\ref{zetataut}) to   
(\ref{transport}) will be called 
the ``natural" solution and the transport
defined by
(\ref{transrho})  the ``natural" transport.

\smallskip

Let $\{\psi_t\,|\,t\in[0,1]\}$ and $\{\tilde\psi_t  \,|\,t\in[0,1]\}$ be two isotopies with
$\psi_1=\tilde\psi_1$. We have ${\cal T}_{\psi}(\rho)$ and
${\cal T}_{\tilde\psi}(\rho)$, 
the sections resulting of the transport of     $\rho$
along both these isotopies. In general
${\cal T}_{\psi}(\rho)$ and
${\cal T}_{\tilde\psi}(\rho)$ will not be equal.
That is, the natural transport is not flat.
In Section 3 we will analyse the corresponding ``curvature".  

\smallskip

 The operator $-D_{X_t}-2\pi if_t$ can be consider from another point of view. 
One can associate to each $C^{\infty}$ function $f$ on $M$ a linear
 operator ${\cal P}_f$ on the space $\Gamma(L)$,
defined by 
$${\cal P}_f(\sigma)=-D_{X_f}\sigma -2\pi if\sigma,$$ 
where $X_f$ is the Hamiltonian vector field determined by $f$.
It is easy to check ${\cal P}_{\{f,g\}}={\cal P}_f\circ {\cal P}_g -{\cal P}_g\circ {\cal P}_f=:[{\cal P}_f,{\cal P}_g]$;  
where the Poisson bracket $\{f,g \}$ is defined as $\omega(X_g,X_f)$. So 
${\cal P}$ is 
  a representation of the Lie algebra 
$C^{\infty}(M)$, the 
{\em prequantization representation}  \cite{jS80}.
  On the other hand, in the algebra of linear operators on $\Gamma(L)$
 one can consider the ideal ${\Bbb C}$ consisting
 of the operators multiplication by a constant, this allows us to define a  representation
 of the algebra
$\text{Lie}(\text{Ham}(M))$   in the algebra 
$\text{End}(\Gamma(L))/{\Bbb C}$.     It is reasonable to conjecture the existence of
  obstructions   to extend the above representation   to
a projective representation 
$$\text{Ham}(M)\rightarrow PL(\Gamma(L))$$
of the group $\text{Ham}(M)$. In Section 3 we will relate these obstructions with the
curvature of the natural transport.

\smallskip

{\it Relation with the Berry phase.}

The connection on the ${\Bbb C}^{\times}$-principal bundle 
$L^{\times}=L-\{ \hbox{zero section}\}$, associated to the prequantum
 bundle $L$, will be denoted by $\alpha$.
Given   $c\in {\Bbb C}$, the vertical vector field on
$L^{\times}$ generated by $c$ will be denoted by $W_c$. That
is, $W_c(q)$ is the vector defined by the curve in $L^{\times}$ given by $q\cdot e^{2\pi ict}$.

Henceforth in this Section  we assume that $F$ is a Lagrangian foliation. 
 Given $\tau\in{\cal Q}_F$,  as $\tau$
 is parallel along the leaves
of the distribution $F$, if $N$ is a leaf of $F$ and    
   if $\tau_{|N}\ne 0$, then $\tau(p)\ne 0$ for all $p\in N$. So 
 $\tau(N)$ is a Planckian submanifold \cite{jS70} of  $L^{\times}$ over $N$.

The proof of the following Lemma is straightforward

\begin{Lem}\label{P:Berry}
If $X\in T_mN$ and $\tau\in{\cal Q}_F$, the vector $\tau_{*}(X)\in T_qL^{\times}$, 
where $q=\tau(m)$,
 satisfies
$\tau_{*}(X)=H(X)(q)+(D_X\tau)(m)$, with $H(X)(q)$ the horizontal lift of $X$ at the point $q$ 
\end{Lem}

\smallskip

Given a Hamiltonian isotopy $\psi_t$ 
and $\tau\in{\cal Q}_F$,
let $\tau_t$ be the family generated by the transport of $\tau$ along $\psi_t$. If  $p\in N$ 
 one can consider in $L^{\times}$  the following curve 
 $$t\rightarrow \tau_t(\psi_t(p))$$
\begin{Prop}\label{P:Berry2}
The tangent vector defined by $\{\tau_t(\psi_t(p))\}_t$ at $q=\tau_u(\psi_u(p))$ is 
$$H(X_u)(q)+W_{-f_u(\pi(q))}(q).$$
\end{Prop}
\begin{pf} For $t$ in a small neighbourhood of $u$, as 
$$\frac{d\tau_t}{dt}=-2\pi if_t\tau_t-D_{X_t}\tau_t,$$
one has
\begin{align}
\tau_t(\psi_t(p))=&\tau_u(\psi_t(p)) \notag \\
&-(t-u)\Big(2\pi if_u(\psi_t(p))\tau_u(\psi_t(p))+(D_{X_u}\tau_u)(\psi_t(p))   \Big)+O((t-u)^2) \notag
\end{align}
This curve defines at $t=u$ the following vector of $T_{q}L^{\times}$
\begin{equation}\label{PRoeq}
(\tau_u)_*(X_u(s))-\big(2\pi if_u(s)\tau_u(s)+(D_{X_u}\tau_u) (s)  \big),
\end{equation}
where $s:=\psi_u(p)$.
As $\tau_u(s)=q$ by the  Lemma \ref{P:Berry}
$(\tau_u)_*(X_u(s))=H(X_u(s))(q)+ (D_{X_u}\tau_u)(s)$. So the expression 
(\ref{PRoeq}) is equal to
$$H(X_u(\pi(q)))(q)-W_{f_u(\pi(q))}(q).$$
In short, the tangent vector at $q$ defined  by the   curve considered is $H(X_u)+W_{-f_u}$.
\end{pf}

We will use the following simple Lemma
\begin{Lem}\label{2trivial}
If $N$  is a connected submanifold of $M$ and
 $\sigma$ and $\rho$ are 
sections of $L$ parallel along $N$, where $\rho$  non identically zero on $N$, 
then   
$\sigma_{|N}=k\rho_{|N}$, with $k$ constant.
\end{Lem}
\begin{pf}
As $\rho$ is parallel along $N$, $\rho(x)\ne 0$ for all  $x\in N$; so there is a function $h$ on $N$ with
$\sigma_{|N}=h\rho_{|N}$.  
 The relation
$$D_A(\sigma_{|N})=A(h)\rho_{|N}+hD_A(\rho_{|N})$$
for every $A\in TN$, implies that $h$ is constant on $N$. 
\end{pf}

Given $\tau\in{\cal Q}_F$, and 
 $\psi_t$  a Hamiltonian {\em closed} isotopy, i. e. such that $\psi_1=\text{id}$,  
 then ${\cal T}_{\psi}(\tau)$ is also $F$-polarized. If $N$ is a leaf of $F$ 
and $\tau_{|N}\ne 0$
 by Lemma \ref{2trivial}
\begin{equation}\label{previousBerry}
{\cal T}_{\psi}(\tau)_{|N}=\kappa\tau_{|N},
\end{equation}
where $\kappa$ is a constant. 
From linearity of the transport and Lemma \ref{2trivial} it follows that 
 $\kappa$ is independent of the section $\tau$.
Hence $\kappa$ can be considered as the {\em holonomy} of  the natural transport,
along the closed isotopy $\psi_t$,
 of $F$-polarized
sections of $L_{|N}$. In Section 3 we will prove the existence of holonomy for the transport 
of {\em arbitrary} sections of $L$.

Now we recall some results of Weinstein about the Berry phase (for details 
see   \cite[page 142]{aW90}). If $\{N_t\}_t$ is a loop of Lagrangian submanifolds
generated by the closed isotopy $\psi_t$.
Let $\epsilon_t$ be a smooth density on $N_t$ such that $\int_{N_t}f_t\epsilon_t=0$.
  Let $\{r_t \}$ be
the family of isomorphisms   
of $(L^{\times}_{|N},\alpha)$ to $(L^{\times}_{|N_t},\alpha)$  determined 
by $\{f_t\}$; that is, 
the   isomorphisms generated by the vector fields
\begin{equation}\label{E:Berryfields}
H(X_t)+W_{-f_t},
\end{equation}
where $H(X_t)$ is the horizontal lift of $X_t$.
The submanifold  $r_1(\tau(N))$ ``differs" from $\tau(N)$ by 
and element $\theta\in U(1)$, that is,
\begin{equation}\label{theta}
r_1(\tau(N))=\theta\tau(N).
\end{equation}
If we denote by $\text{hol}$ the holonomy on $N$ defined 
by de connection $\alpha$, 
$\,\text{hol}:\pi_1(N)\rightarrow U(1)$, then the Berry 
phase of the family $(N_t,\epsilon_t)$ of weighted submanifolds
is the class of $\theta$ in the quotient 
$U(1)/(\hbox{Im}(\text{hol}))$. Up to here the results of Weinstein.

\smallskip

\begin{Thm}\label{TheBerph}
If $\psi_t$ is a closed Hamiltonian isotopy and $N$ a connected leaf of the Lagrangian 
foliation $F$, then 
the Berry phase of $(N_t,\epsilon_t)$, with $N_t=\psi_t(N)$,  is the class in 
$U(1)/(\hbox{Im}(\text{hol}))$ of the holonomy of  the natural  transport
along $\psi_t$ 
of $F$-polarized sections of $L_{|N}$. 
 \end{Thm}
\begin{pf}
Given $p\in N$, by  Proposition \ref{P:Berry2} and (\ref{E:Berryfields})
 the curves in $L^{\times}$ 
$\{\tau_t(\psi_t(p) \}_t$ and $\{r_t(\tau(p)) \}_t$   define the same vector field. As they
take the same value for $t=0$, it turns out that
 $r_t(\tau(p))=\tau_t(\psi_t(p))$, for all $p\in N$; hence the above complex $\theta$ in 
(\ref{theta}) is determined by
\begin{equation}\label{trivial}
\tau_1(p)=\theta\tau(p).
\end{equation}
From (\ref{previousBerry}) and (\ref{trivial})
we conclude $\theta=\kappa$.
\end{pf}

\section{The holonomy of the natural transport}

We will prove that it makes sense to define the holonomy of
the natural transport of {\em arbitrary} sections of $L$.
We  start with  $\psi=\{\psi_t\,| \, t\in[0,1] \}$ a {\em closed} Hamiltonian isotopy, generated by the
time dependent Hamiltonian function $f_t$;  
 that is, $\psi$ is a loop at $\text{id}$ in the group
$\text{Ham}(M)$. So we must consider the corresponding equation 
(\ref{transrho}) and study its solution. Let $\mu$ be
a local frame for $L$, defined on $R\subset M$, and $\beta$ the connection form in this frame. 
There is a time dependent function $g(t,\,.\,)$ such that
\begin{equation}\label{grho}
\rho_t(p)=g(t,p)\mu(p),\,\,\, t\in[0,1],\,\,\,\text{and}\,\,\, p\in R.
\end{equation}
Hence   (\ref{transrho})
can be written
\begin{equation}\label{basreduc}
\frac{\partial}{\partial t}g(t,\,.)=
-X_t\big(g(t,\,.\,)\big)-\beta(X_t)g(t,\,.\,)-2\pi if_{t} g(t,\,.\,),\,\,\,\, g(0,p)=p \,\,\text{for all}\,\, p\in R.
\end{equation}

Fixed a point $q\in M$  we put  $\sigma(t):=\psi_t(q)$.
Assumed that the closed curve $\sigma$  is contained
 in $R$, the equation (\ref{basreduc}) on the points of this curve
is
\begin{equation}\label{basreduc1}
\frac{\partial g}{\partial t}(t,\sigma(t))+X_t(\sigma(t))\big(g(t,\,.\,)\big)=
-\beta_{\sigma(t)}(X_t)g(t,\sigma(t))-2\pi if_{t}(\sigma(t))g(t,\sigma(t))
\end{equation}
The second sumand on the left hand side is the action of the vector $X_t(\sigma(t))$ on
the function $g(t,\,.\,):M\rightarrow {\Bbb R}$.
If we consider the curve $\hat \sigma:[0,1]\rightarrow {\Bbb R}\times M$, defined by
$\hat\sigma(t)=(t,\sigma(t))$ and we put $\hat g(t):=g(\hat\sigma(t))$, equation (\ref{basreduc1})
can be written
\begin{equation}
\frac{d\,\hat g}{dt}=-\beta_{\sigma(t)}(X_t)\hat g(t)-2\pi if_{t}(\sigma(t))\hat g(t)
\end{equation}
Hence
$$\hat g(t)=g(0,q)\,\text{exp}\Big(\int_0^t\big(-\beta_{\sigma(t')}(X_{t'})-
2\pi if_{t'}(\sigma(t'))\big)dt'\Big).$$
 As the closed curve $\sigma$ is nullhomologous  \cite[page 334]{M-S98}, let $S$ be any oriented $2$-submanifold bounded
by the closed curve $\sigma$, then   
 $$\int_0^1\beta_{\sigma(t)}(X_t)\,dt=\int_S d\beta.$$
As the curvature of $L$
is $-2\pi i\omega$ we have
   \begin{equation}\label{gsombrero}
\hat g(1)= g(0,q)\,\text{exp}\Big(2\pi i\int_S\omega-2\pi i\int_0^1f_{t}(\psi_t(q))\,dt   \Big).
\end{equation}

Given the loop $\psi$ in $\text{Ham}(M)$, the Hamiltonian vector fields $X_t$ determine the 
Hamiltonian $f_t$ up to an additive constant. In certain cases it is possible to fix this Hamiltonian
function in a natural way; for instance when $X_t$ is an invariant
 vector field on a coadjoint orbit.
 In a general case $f_t$ can be fixed by imposing
the condition that $f_t$ has zero mean with respect to  the canonical measure
on $M$ induced by $\omega$; henceforth we assume that $f_t$ satisfies this
{\em normalization condition}.

For $p$ point in $M$ one defines
\begin{equation}\label{kappa}
\kappa_p(\psi):=\text{exp}\Big(2\pi i\int_S\omega-2\pi i\int_0^1f_{t}(\psi_t(p))\, dt   \Big),
\end{equation}
where   $S$ is any surface bounded by the closed curve $\psi_t(p)$ in $M$.

Given the closed isotopy $\psi$, one can define 
the {\em action integral} \cite{aW89} \cite{M-S98} 
${\cal A}(\psi)(p)$ around the curve $\psi_t(p)$ as the element of ${\Bbb R}/{\Bbb Z}$ determined
by $(2\pi i)^{-1}$ times the exponent of (\ref{kappa}). 
Hence $\kappa_p(\psi)=\text{exp}(2\pi i{\cal A}(\psi)(p))$.

If the Hamiltonian function is independent of $t$ (i.e., the loop $\psi$ is $1$-parameter
subgroup in $\text{Ham}(M)$), then it is constant along
$\psi_t(p)$.  Consequently the  second integral in (\ref{kappa}) is equal to $f(p)$.

From (\ref{grho}) it follows $\rho_1(q)=\kappa_q(\psi)\rho(q)$. And
by choosing appropriate local frames, one can prove for any $\rho\in\Gamma(L)$ 
 \begin{equation}\label{tau1=}
\rho_1(p)=\kappa_p(\psi)\rho(p),\,\,\,\text{for any}\,\, p\in M
\end{equation}

To study the function $\kappa(\psi):M\rightarrow U(1)$ we will use properties 
 of the sections of $L$ polarized with respect
to certain foliations. There may be topological obstructions to the existence
of such foliations; however we will prove some properties of that function using
the existence of families of vector fields which define foliations in  parts of $M$.
  
Let ${\cal B}:=\{B_1,\dots,B_m\}$ be a set of vector fields on $M$ which define an 
$m$-dimensional foliation  on $M-K$, where $K$ is a subset of $M$.  
 This foliation will be denoted also by ${\cal B}$. We put
 $${\cal B}_t=\{B_j(t):=(\psi_t)_*(B_j) \}_{j=1,\dots,m},$$ and 
this set defines    a foliation  
on $M-\psi_t(K)$. Moreover if $N$ is a leaf of ${\cal B}$, then $N_t=\psi_t(N)$ is a 
leaf of ${\cal B}_t$.
On the other hand, according to Proposition \ref{ProposiNueva}, if $\tau$ is a section of $L$ which is ${\cal B}$-polarized, 
that is, such that 
$D_{B_j}\tau=0$,  $\,j=1,\dots,m$,
then the section $\tau_t$  solution to (\ref{transrho}), 
  is ${\cal B}_t$-polarized.

Let $N\subset M-K$ be a {\em connected} integral submanifold of ${\cal B}$. 
Given $\tau$ a ${\cal B}$-polarized
 section of $L_{|N}$  with $\tau$ non identically zero on $N$,
then $\tau_1(p)=\kappa_p(\psi)\tau(p)$ for $p\in N$. As $\tau_1$ and $\tau$ are
${\cal B}$-polarized by Lemma \ref{2trivial}
one deduces that $\kappa_p(\psi)$ is independent of the point $p\in N$.
The above results can be summarised in the following 
\begin{Prop}\label{teoremaindepe} Let  $(M,\omega)$ be a compact,  
 quantizable manifold, and $\psi$ a loop in $\text{Ham}(M)$ at $\text{id}$. If $p,q$  
are points
which belong to a connected integral submanifold $N$ of the foliation on $M-K$ defined
by ${\cal B}$, then
 $\kappa_q(\psi)=\kappa_p(\psi)$, provided that there is a ${\cal B}$-polarized   section of $L$ 
 non-zero on $N$. 
 \end{Prop}

\begin{Cor}\label{CoroThe}  
 If $N^i$ $\, (i=1,2)$  is a connected integral submanifold of ${\cal B}^i$,
and $\tau^i$ a ${\cal B}^i$-polarized section of $L_{|N^i}$ with $\tau^i\ne 0$. Then 
$\kappa_{q_1}= \kappa_{q_2}$, if $q_i\in N^i$ and $N^1\cap N^2\ne\emptyset$.
\end{Cor}
\begin{pf} If $p\in N^1\cap N^2$, then for any loop $\psi$ one has
$\kappa_{q_1}(\psi)= \kappa_{p}(\psi)= \kappa_{q_2}(\psi)$. 
\end{pf}

On the other hand, if $N$ is a simply connected integral submanifold 
of an isotropic foliation ${\cal B}$, as $\omega_{|TN}=0$, the parallel transport determined 
by the connection of $L$ allows us to define a nonzero section $\rho$ of $L_{|N}$ 
parallel along $N$.   This fact permits other formulation of
 Proposition \ref{teoremaindepe} without assuming the existence of the nonzero 
polarized section.
\begin{Prop}Let us suppose that $(M,\omega)$ 
 is a compact,  
 quantizable manifold, 
and that
 $p$ and $q$ are two points which belong to
a connected integral submanifold $N$ of the isotropic foliation ${\cal B}$,
   if $N$ can be written as a finite union
of open simply connected  subsets, then $\kappa_{p}=\kappa_{q}$. 
\end{Prop}
\begin{pf} It is a consequence of the preceding remark and Corollary \ref{CoroThe}.
\end{pf}

Corollary \ref{CoroThe} admits also a similar version without supposing the
existence  of $\tau^i$, if we assume that $N^i$ can be expressed as a finite union of 
simply connected open subsets. So we have
\begin{Prop}\label{Thmfinal}
Assumed that $(M,\omega)$ is a compact and quantizable manifold. 
Let ${\cal B}^i$ (i=1,2) be two sets of vector fields which define
isotropic foliations on $M-K^i$,  and $N^i$ a connected integral
 submanifold of ${\cal B}^i$,
if $N^i$ can be written as a finite union
of simply connected  open subsets and $N^1\cap N^2\ne\emptyset$,
then 
$\kappa_{q_1}= \kappa_{q_2}$, for $q_i\in N^i$. 
\end{Prop}

 \smallskip 

Let $Y$ be a transversal vector field on $M$, that is, $Y$ is a section of $TM$
which is transversal to the zero section of $TM$. Then the Euler class 
$e(M)\in H^{2n}(M)$ of $M$ is 
Poincar\'e dual of the zero locus of $Y$; so this zero locus is a finite set $K$
of points of $M$. And from the transversality theory we conclude that this property
 is also valid for any ``generic" vector field. Each point of  $M-K$ belongs to 
a non constant integral curve of $Y$. If $p$ and $q$ are two arbitrary points in $M$ one
 can choose generic vector fields 
 $Y_1,\dots, Y_{m}$ on $M$ such that $p$ and $q$ can be joined by a path which is the 
juxtaposition of   curves, each of which is a non constant integral curve of 
some $Y_j$. As the curves are isotropic submanifolds of $M$
 by Proposition \ref{Thmfinal},
 $\kappa_p=\kappa_q$.

If $\xi$ and $\psi$ are two loops in $\text{Ham}(M)$ based at $\text{id}$ we can define 
$\xi\cdot\psi$ as the loop given by the usual product of paths, it is immediate to check that
$\kappa(\xi\cdot\psi)=\kappa(\xi)\kappa(\psi)$. 

 By (\ref{tau1=}) and the foregoing reasoning  we can state the following
\begin{Thm}\label{Therefi} If $(M,\omega)$ is compact and quantizable,   
the correspondence 
$$\kappa:\{\text{Loops in Ham(M) based at id}  \}\rightarrow U(1)$$
 defined by
$$\kappa(\psi)=\text{exp}\Big(2\pi i\int_S\omega-2\pi i\int_0^1f_{t}(\psi_t(q))dt   \Big),$$
 $q$  being an arbitrary point of $M$ and $S$ any surface bounded by the closed curve
$\{\psi_t(q)\}$, is a well-defined map which satisfies
$\kappa(\xi\cdot\psi)=\kappa(\xi)\kappa(\psi)$. Moreover 
\begin{equation}\label{calTkappa}
{\cal T}_{\psi}\rho=\kappa(\psi)\rho,
\end{equation}
for any section $\rho$ of the prequantum bundle $L$.
\end{Thm}

By (\ref{calTkappa}) it makes sense to call $\kappa(\psi)$ the {\em holonomy}
of the natural transport along the loop $\psi$.

\begin{Cor} The action integral ${\cal A}(\psi)(p)$ is independent of $p$.
\end{Cor}

\begin{Cor}\label{Weinsteinresult}  
Let $f$ be a Hamiltonian function  such that it defines a 
$1$-parameter loop $\{\psi_t \,|\,t\in[0,1] \}$
of symplectomorphisms; if $p$ is critical point of $f$, then 
 $\kappa(\psi)=\text{exp}(-2\pi if(p))$.
If $p$ and $q$ are critical points of $f$ then 
$f(p)=f(q)\,\,(\text{mod}\,{\Bbb Z})$.
\end{Cor}
\begin{pf} As $\psi_t(p)=p$ for all $t$, the corollary is a consequence of (\ref{kappa}). 
\end{pf}
\smallskip
This relation among the critical values of $f$ has been proved in \cite{aW89} using
the invariant ${\bf A}$ mentioned in the Introduction.

The relation (\ref{calTkappa})   and Theorem
\ref{TheBerph} imply
\begin{Cor} Let $\{N_t:=\psi_t(N)\}$, with $N$ a connected, simply connected Lagrangian 
submanifold of $M$ and $\psi=\{\psi_t\}$ a loop in $\text{Ham}(M)$,
then the Berry phase of the family  $(N_t,\epsilon_t)$ of weighted submanifolds is $\kappa(\psi)$. 
\end{Cor}

\smallskip

Next we   will study the behaviour of   $\kappa(\psi)$ 
under $C^{1}$-deformations of $\psi$. Let $\psi=\{ \psi_t\,|\, t\in[0,1]\}$ be
a loop in $\text{Ham}(M)$ with $\psi_0=\psi_1=\text{id}.$ We consider
the derivative of $\kappa(\psi^s)$ with respect to the parameter $s$ in a deformation
$\psi^s$ of $\psi$. That is, $\psi^s=\{\psi^s_t\,|\, t\in[0,1] \}$ is an isotopy with
$\psi^s_0=\psi^s_1=\text{id}$ generated by the time dependent Hamiltonian $f^s_t$; 
furthermore we assume that
  $\psi^0=\psi$. By $\{X^s_t\}_t$ is denoted the family of Hamiltonian vector fields 
defined by $\{f^s_t\}_t$.

For $q\in M$ we put $\sigma^s(t):=\psi^s_t(q)$, so
 $\{\sigma^s(t)\,|\, t\in[0,1]\}$ is a closed curve and then
$$\kappa(\psi^s)=\text{exp}\,\Big(2\pi i\int_{S^s}\omega
  -2\pi i\int_0^1f_t^s(\sigma^s(t))\Big)=:\text{exp}(2\pi i\Delta(s)),$$              
where $S^s$ is a surface bounded by the curve $\sigma^s$. 
We set 
$$X_t:=X^0_t,\,\,\,  f_t:=f^0_t,\,\,\,  \sigma(t)=\sigma^0(t) .$$
The variation of $\sigma^s(t)$ with $s$
permits to define the vector fields $Y_t$; that is,
\begin{equation}\label{Yts}
Y_t(\sigma^s(t)):=\frac{\partial}{\partial s} \sigma^s(t).            
\end{equation}
For an ``infinitesimal" $s$ the curves $\sigma^l$, with $l\in[0,s]$ determine the 
``lateral surface" $J$ of one ``wedge"  whose base and cover
 are the surfaces $S$ and $S^s$ respectively. The ordered pairs of vectors   
$(X_t(\sigma(t)),Y_t(\sigma(t)))$
fix an orientation on $J$, which   in turn 
determines an orientation on the closed surface   $T=S\cup J\cup S^s$.
If we assume that $S$ and $S^s$ are oriented by means of the orientations 
of curves $\sigma$ and $\sigma^s$, from the fixed orientation on $T$ it 
 follows $T=J-S+S^s$. 

As $\omega$ satisfies the integrality condition
 \begin{equation}\label{derivadaS}
-\int_S\omega+\int_{S^s}\omega  =-\int_{J}\omega \,\,\text{(modulo}\,{\Bbb Z}\text{)}.
\end{equation}
Moreover 
\begin{equation}\label{intJOmega}
\int_{J}\omega=
s\int_0^1\omega\big(X_t(\sigma(t)),\,Y_t(\sigma(t)) \big)dt +O(s^2).
\end{equation}
On the other hand, for a given   $t\in[0,1]$  
\begin{equation}\label{derivfst}
\Big(\frac{d}{ds} f^s_t(\sigma^s(t))  \Big)_{|s=0}  =
\Big(\frac{\partial}{\partial s}f^s_t(\sigma(t))\Big)_{|s=0}+Y_t(\sigma(t))(f_t).
\end{equation}
We set
$$\Dot f_t(p):=\Big(\frac{\partial}{\partial s}f^s_t(p)\Big)_{|s=0}.$$
 As $\iota_{X_t}\omega=-df_t$, from
(\ref{derivfst}) it follows
\begin{equation}\label{d|dsint}
\frac{d}{ds}\Big|_{s=0} \int_0^1f^s_t(\sigma^s(t))dt = \int_0^1\Dot f_t(\sigma(t)) dt
-\int_0^1\omega\big(X_t(\sigma(t)),Y_t(\sigma(t))\big)dt.
\end{equation}
 By (\ref{derivadaS}), (\ref{intJOmega}) and  (\ref{d|dsint})
$$\Delta(s)-\Delta(0)=-s\int_0^1\Dot f_t(\sigma(t)) dt +O(s^2)\,\, \text{(modulo}\,{\Bbb Z}\text{)}.$$
So 
$$\kappa(\psi^s)-\kappa(\psi)=-2\pi is\kappa(\psi)\int_0^1\Dot f_t(\sigma(t)) dt+O(s^2),$$
 and finally
\begin{equation}\label{derivkappa}
\Big(\frac{d}{ds}\kappa(\psi^s)\Big)_{|s=0}=-2\pi i \kappa(\psi)\int_0^1\Dot f_t(\psi_t(q)) dt.
\end{equation}

\smallskip

 By ${\cal L}$ is denoted the space of $C^1$-loops in $\text{Ham}(M)$
based at $\text{id}$; that is, ${\cal L}$ is the space of isotopies ending at $\text{id}$. Given $\psi\in {\cal L}$, let $\psi^s$ a curve in 
${\cal L}$ with $\psi^0=\psi$. For each $s$ one has the corresponding time dependent
Hamiltonian function $f^s_t$. The tangent vector $Z$ defined by $\psi^s$ is determined
by the family of functions 
$$\Dot f_t:=\frac{\partial}{\partial s}\biggr|_{s=0}f^s_t,$$
which in turn can be identified with the corresponding Hamiltoninan family 
 of vector fields.
 
On ${\cal L}$ we define the $1$-form $\Omega$ as follows: Given $Z\in T_{\psi}{\cal L}$,
determined by the family $\{\Dot f_t\},$
\begin{equation}\label{defiOMEGA}
\Omega_{\psi}(Z):=\int_0^1\Dot f_t(\psi_t(q))dt,
\end{equation}
where $q$ is any point of $M$. The left hand side in (\ref{derivkappa}) 
is independent of the point $q$, and so the right
hand side is also; therefore $\Omega$ is well-defined.

 If  $\Omega=0$, then for any loop $\psi$ in $\text{Ham}(M)$ and any
deformation $\psi^s$ of $\psi$ we have 
$$\Big(\frac{d}{ds}\kappa(\psi^s)  \Big)_{|s=0}=0,$$
and conversely. In this case
  $\kappa$ is invariant under homotopies.

The Lie algebra of the group $\text{Ham}(M)$ consists of all smooth functions on $M$
which satisfy the normalization condition. The prequantization  map ${\cal P}$ is a representation
of this algebra, as we said in Section 2.  In general ${\cal P}$ is
 not the tangent representation
of one representation of $\tilde{\text{Ham}}(M)$, the universal cover of $\text{Ham}(M)$. In the following we analyse
this issue.
An element of $\tilde{\text{Ham}}(M)$ is a homotopy class of a curve in $\text{Ham}(M)$ which
starts at $\text{id}$, i. e. the homotopy class $[\psi]$ of a Hamiltonian isotopy $\psi$.
When $\Omega$ vanishes ${\cal T}_{\psi}$ depends only on the homotopy class $[\psi]$, this fact
allows to construct a representation of $\tilde{\text{Ham}}(M)$ whose 
tangent representation is ${\cal P}$.

\begin{Prop} If $\Omega=0$, then the prequantization map ${\cal P}$ extends to a representation of
$\tilde{\text{Ham}}(M)$.
\end{Prop}
\begin{pf}  Given the isotopy $\psi$, let   $\{ \psi^s\}$ be a
 deformation of $\psi$. For each $s$ the path $\zeta^s$ in $\text{Ham}(M)$ defined as the usual product path
$\psi^s\cdot\psi^{-1}$ of the corresponding paths is a closed isotopy. Since
  $\Omega=0$,  
$$\Big(\frac{d}{ds} \kappa(\zeta^s)\Big)_{|s=0}=0.$$
As 
$$\big( {\cal T}_{\psi^{-1}}\circ {\cal T}_{\psi^{s}}  \big)(\rho)=
{\cal T}_{\psi^s\cdot\psi^{-1}}(\rho)=\kappa(\psi^s\cdot\psi^{-1})\rho,$$
for every $\rho\in\Gamma(L)$, then  
$$\frac{d}{ds}\biggr|_{s=0}\big( {\cal T}_{\psi^{-1}} {\cal T}_{\psi^{s}}  \big)(\rho)=
\Big(\frac{d}{ds}\kappa(\zeta^s)\Big)_{|s=0}\rho=0.$$
So the transport ${\cal T}_{\psi}$ along $\psi$ 
depends only on the homotopy class $[\psi]$.
That is,  ${\cal T}$
is well-defined on $\tilde{\text{Ham}}(M)$. 

If  $\psi$ and $\chi$ are isotopies, one can consider $\chi\circ\psi$,
the isotopy defined by $(\chi\circ\psi)_t=\chi_t\circ\psi_t$. On the other hand one has
the juxtaposition $\psi\star\chi$ given by $(\psi\star\chi)_t=\psi_{2t}$, for $t\in[0,0.5]$
and $(\psi\star\chi)_t=\chi_{2t-1}\circ\psi_1$, for $t\in[0.5,1]$. As
$[\chi\circ\psi]=[\psi\star\chi]$ (see \cite{M-S98}) we have
$${\cal T}_{[\chi][\psi]}={\cal T}_{[\psi\star\chi]}={\cal T}_{[\chi]}\circ {\cal T}_{[\psi]}.$$
Hence ${\cal T}$ is a representation of $\tilde{\text{Ham}}(M)$ and, by construction,
its tangent representation is ${\cal P}$.
 \end{pf}

In a similar way one can prove the following 
\begin{Thm} Let $(M,\omega)$ be a compact, quantizable manifold. The following
properties are equivalent

(i) The $1$-form $\Omega$ vanishes.  

(ii) For any simply connected Lagrangian submanifold $N$ of $M$ and every loop ${\psi_t}$ in $\text{Ham}(M)$
at $\text{id}$ the Berry phase of the loop $\{\psi_t(N) \}_t$ of Lagrangian 
submanifolds depends only  on the homotopy class of  ${\psi_t}$.

(iii) Given an arbitrary foliation $F$ of $M$ and an arbitrary Hamiltonian isotopy $\psi$, the
natural identifications of ${\cal Q}_F$ and ${\cal Q}_{\psi_1(F)}$ defined by
${\cal T}_{\psi}$ and  ${\cal T}_{\psi'}$ are equal, for all $\psi'\in[\psi]$.

\end{Thm}

\begin{pf} We assume (i). If $\psi^s$ is a deformation of $\psi\in{\cal L}$, as in the
 foregoing Proposition $\psi^{s}\cdot\psi^{-1}$ is a closed isotopy. By (i) 
$$\Big(\frac{d\, \kappa(\psi^s\cdot\psi^{-1})}{ds}\Big)_{|s=0}=0.$$
 From (\ref{calTkappa}) and Theorem \ref{TheBerph} it follows property (ii).  

Conversely, let $\psi$ be an element of ${\cal L}$ and  $\psi^s$ an 
arbitrary curve in ${\cal L}$
with $\psi^0=\psi$. This curve defines a deformation of $\psi$. Let us take $\tau\in{\cal Q}_F$
with $\tau_{|N}\ne 0$, for $N$ a leaf of a Lagrangian foliation $F$. By (ii) and Theorem \ref{TheBerph}
 $\kappa(\psi)\tau_{|N}= \kappa(\psi^s)\tau_{|N}$ for all $s$. Therefore
 $\Big(\frac{d\,\kappa(\psi^s)}{ds}\Big)_{|s=0}=0$; consequently
$\Omega_{\psi}=0$. 
\end{pf}
 
\smallskip

 Next we study a particular case: The behaviour of $\kappa(\psi)$
under    deformations consisting of $1$-parameter subgroups. Let us suppose
that $\psi^s$ for each $s$ is a 
$1$-periodic Hamiltonian flow;
 then $f^s_t$ is independent of $t$ and we put $f^s_t=f^s.$ One defines
 the function $\Dot f$ by $\Dot f(p)=\Big( \frac{d}{ds}f^s(p)\Big)_{|s=0}.$
As $\{\sigma^s(t)   \,|\, t\in[0,1] \}$ is an integral curve for the Hamiltonian function
$f^s$
$$f^s(\sigma^s(t))=f^s(q)=(f+s\Dot f)(q)+O(s^2)=f(q)+s\Dot f(q)+O(s^2).$$
On the other hand
$$f^s(\sigma^s(t))=(f+s\Dot f)(\sigma^s(t))+O(s^2)=f(q)+s\big(Y_t(\sigma(t))(f)+\Dot f(\sigma(t))\big)+O(s^2).$$
Therefore
$$\Dot f(q)=Y_t(\sigma(t))(f)+\Dot f(\sigma(t)).$$
Now $df=-\iota_X\omega$, then 
$$\int_0^1\Dot f(\sigma(t))dt=\Dot f(q)+\int_0^1\omega\big(X(\sigma(t)),Y_t(\sigma(t))\big)dt.$$
The symplectomorphism $\delta:=\psi^s_t\circ\psi_t^{-1}$ applies the curve $\sigma(t)$ into
$\sigma^s(t)$. Hence $\delta(S)$ is a surface whose boundary is $\sigma^s(t)$ and
\begin{equation}\label{punto}
\int_{S^s}\omega=\int_S\delta^*\omega=\int_S\omega.
\end{equation}
From (\ref{punto}), (\ref{derivadaS}) and (\ref{intJOmega}) it follows
$$\int_0^1\omega\big(X(\sigma(t)),Y_t(\sigma(t))\big)dt=0.$$
From (\ref{derivkappa}) it follows
\begin{equation}\label{deri2}\frac{-1}{2\pi i \kappa(\psi)}
\Big(\frac{d}{ds}\kappa(\psi^s)\Big)_{|s=0}=
\Dot f(q).
\end{equation}
As the left hand side in (\ref{deri2}) is independent of the point $q$, it 
turns out that $\Dot f$ is constant on $M$. The normalization condition of each $f^s$
implies
$$0=\int_M f^s\omega^n=\int_M(f+s\Dot f)\omega^n
+O(s^2)=s\Dot f\int_M\omega^n+O(s^2).$$
Hence $\Dot f\equiv 0$, and by (\ref{derivkappa})
$$\Big(\frac{d}{ds}\kappa(\psi^s)  \Big)_{|s=0}=0.$$
One has
\begin{Thm}\label{invakappahomot} $\kappa$ is invariant under homotopies consisting of $1$-parameter subgroups in $M$.
\end{Thm}

\begin{Cor} Let $\psi$ and $\psi'$ be $1$-periodic Hamiltonian flows 
  generated by the Hamiltonian functions
$f$ and $f'$ respectively. If $\psi$ and $\psi'$ are homotopic in the space of $1$-parameter
subgroups, then 
$$ f(p)=f'(p')\,\,\text{(mod}\,{\Bbb Z}\text{)},$$
for $p$ and $p'$ critical points of $f$ and $f'$ respectively.
\end{Cor}
\begin{pf} It is a consequence of Theorem \ref{invakappahomot} and 
Corollary \ref{Weinsteinresult}.  
\end{pf}

 \section{A grading in $\pi_2(\text{Ham}(M))$.}

We will prove in this Section that the $1$-form $\Omega$ on ${\cal L}$ is closed.
If  $\phi:=\{\phi^s\}$
is a closed curve in ${\cal L}$, one can consider the map
$\kappa(\phi^-):s\in S^1\mapsto\kappa(\phi^s)\in U(1)$; its winding number is 
$$\text{deg}(\kappa(\phi^-))=\int_{S^1}\frac{1}{2\pi i\kappa (\phi^s)}\frac{d\,\kappa(\phi^s)}{ds}\,ds.$$
By (\ref{derivkappa}) and (\ref{defiOMEGA}) this winding number is equal to
$$-\int_{S^1}\Omega_{\phi^s}(\Dot \phi^s)ds,$$
where $\Dot \phi^s$ is the vector of $T_{\phi^s}{\cal L}$ defined
 by the curve $\{\phi^s\}_s$.

If $\phi$ and $\xi$ are two homotopic loops in ${\cal L}$, then there is
a homotopy $_r\phi^s$ such that $_{0}\phi^s=\phi^s$ and  $_{1}\phi^s=\xi^s$.
Therefore $\kappa(_{r}\phi^-)$ is a homotopy between the maps $\kappa(\phi^-)$ and 
$\kappa(\xi^-)$, so these maps have the same degree (see  \cite[page 129]{mH88}).

If $\phi$ and $\varphi$ are loops in ${\cal L}$ based at the same point and 
  $\zeta=\phi\cdot\varphi$ is the path product, then $\text{deg}(\kappa(\zeta^-))$ is equal to
 $$\frac{1}{2\pi i}\Big( \int_0^{0.5}\frac{1}{\kappa(\phi^{2s})}
\frac{d\,\kappa(\phi^{2s})}{ds} ds +
  \int_{0.5}^{1}\frac{1}{\kappa(\varphi^{2s-1})}
\frac{d\,\kappa(\varphi^{2s-1})}{ds} ds 
  \Big),$$
 and this expression   is equal to 
$\text{deg}(\kappa(\phi^-))+ \text{deg}(\kappa(\varphi^-))$.
Thus we have 

\begin{Thm}\label{Omegadegree} $\Omega$ defines an element of $H^1({\cal L},{\Bbb Z})$.
Moreover, if $\phi$ is a closed
curve on ${\cal L}$ then $-\Omega([\phi])$ is the degree of the map   
$\kappa(\phi^-)$.
\end{Thm}

We denote by $c$ the loop in $\text{Ham}(M)$ defined by $c(s)=\text{id}$, for all $s$.
Since $\pi_1({\cal L},\,c)=\pi_2(\text{Ham}(M),\,\text{id})$, 
the form $\Omega$ defines a {\em degree} on $\pi_2(\text{Ham}(M),\text{id})$:
 Given
$[\phi]\in\pi_2(\text{Ham}(M),\,\text{id})$
\begin{equation}
\text{Deg}([\phi]):=\Omega([\phi])=-\text{deg}(\kappa(\phi^-)).
\end{equation}
As $\text{Deg}$ is a homomorphism, this  grading on $\pi_2(\text{Ham}(M))$ is
compatible with the group structure.

If $\Omega$ is exact, then $\text{Deg}=0$. In this case there is a potential map $H:{\cal L}\rightarrow {\Bbb R}$
such that, if $\{ \nu^s\}_s$ is a curve in ${\cal L}$ starting at $c\in{\cal L}$
$$H(\nu^s)= \int_0^s\Omega_{\nu^a}(\Dot\nu^a)\,da=
\frac{-1}{2\pi i}\int_0^s\frac{1}{\kappa(\nu^a)}\frac{d\, \kappa(\nu^a)}{da}\,da.$$
So 
$$\frac{d\,H(\nu^s)}{ds}=\frac{-1}{2\pi i\kappa(\nu^s)}\frac{d\, \kappa(\nu^s)}{ds}.$$
By (\ref{kappa}) $\kappa(c)=1$, so $\kappa(\nu^s)=\text{exp}(-2\pi i H(\nu^s)).$
Hence for 
 every $\psi\in {\cal L}$ that can be joined with $c$ by a path, we have   
$\kappa(\psi)=\text{exp}(-2\pi iH(\psi)).$ A similar expression holds in each
connected component of ${\cal L}$.
Thus $H$ is a lifting of the action integral function 
${\cal A}:{\cal L}\rightarrow {\Bbb R}/{\Bbb Z}$ to an ${\Bbb R}$-valued function.

Conversely, if there is a lifting of ${\cal A}$ to an ${\Bbb R}$-valued function,
then $\text{Deg}=0$; i.e. $\Omega$ is exact. In short
\begin{Prop}
 The class $[\Omega]\in H^1({\cal L},{\Bbb Z})$
is the obstruction to existence of a lifting of ${\cal A}$ to an 
${\Bbb R}$-valued function.
\end{Prop} 

 \smallskip

A generic element of $\pi_2(\text{Ham}(M),\, \text{id})$ is given by a map $\phi=(\phi^s_t)$
from $I^2$ into $\text{Ham}(M)$, such that for each $s$ $\phi^s=\{\phi^s_t \}_t$
 is a Hamiltonian 
isotopy ending at $\text{id}$, defined by the normalized time dependent Hamiltonian $f^s_t$.
One can also consider a family of  particular elements in $\pi_2(\text{Ham}(M),\, \text{id})$,
those $\chi$ such that for  each $s$ $\chi^s$ is the Hamiltonian flow associated to 
a Hamiltonian function.  One has the following result
\begin{Prop}If $[\phi]=[\chi]\in \pi_2(\text{Ham}(M),\, \text{id})$, then
$$\int_0^1\int_0^1\Big(\frac{\partial \, f^s_t}{\partial s}  \Big)\big( \phi^s_t(q)\big) dt\,ds=0,$$
for every $q\in M$.
\end{Prop} 
\begin{pf}
\begin{equation}\label{Omega([chi}
\Omega([\chi])=\Omega([\phi])=\int_0^1\Omega_{\phi^s}(\Dot\phi^s)\,ds
\end{equation}
By (\ref{defiOMEGA})
\begin{equation}\label{Omegaphis(}
\Omega_{\phi^s}(\Dot\phi^s)=\int_0^1
\Big(\frac{\partial \, f^s_t}{\partial s}  \Big)\big( \phi^s_t(q)\big) dt,
\end{equation}
for every $q\in M$.

On the other hand, $\kappa(\chi^s)$ is independent of $s$ by Theorem \ref{invakappahomot}.
So the map $\kappa(\chi^-)$  has degree $0$. The Proposition follows from Theorem \ref{Omegadegree}, (\ref{Omega([chi})
and (\ref{Omegaphis(}).
\end{pf}

\section{Example: Coadjoint orbits of $SU(2)$}

We will check the above results when $M$ is a coadjoint
 orbit \cite{aK76} of the group $SU(2)$.

Let $\eta$ be the element of ${\frak su}(2)^*$
$$\eta:\begin{pmatrix} ai & w \\ -\bar w & -ai 
\end{pmatrix}\in{\frak su}(2)\rightarrow
 ka\in{\Bbb R},$$
 where $k$ is a non-zero real number. It is straightforward to see that the subgroup of isotropy
$G_{\eta}$ of $\eta$ is the subgroup $U(1)$ of $SU(2)$. So the coadjoint
 orbit ${\cal O}_{\eta}$ of $\eta$ can
be identified with $SU(2)/U(1)=S^2$. If $\mu\in{\cal O}_{\eta}$ then $\mu =g\cdot\eta$, with
\begin{equation}\label{xysu2}
g=\begin{pmatrix} x & y \\ -\bar y & \bar x \end{pmatrix}\in SU(2).
\end{equation}
If we put 
\begin{equation}\label{xyinS2}
x=\cos (\theta/2)\,\text{exp}(i\phi_1),\,\,  
y=\sin (\theta/2)\,\text{exp}(-i\phi_2), \,\, \text{with} \,\,0\leq \theta\leq \pi,
\end{equation} 
then
the    point in $S^2$ corresponding to $\mu\in{\cal O}_{\eta}$ through the
diffeomorphism   ${\cal O}_{\eta}\simeq SU(2)/U(1)\simeq S^2$ has the spherical coordinates 
$(\theta,\phi=\phi_1-\phi_2).$ 

On the other hand ${\frak su}(2)=  {\Bbb R}A\oplus{\Bbb R}B \oplus{\Bbb R}Z$, with
$$A=\begin{pmatrix} 0 & i \\i & 0 \end{pmatrix},\,\,
B=\begin{pmatrix} 0 & 1 \\-1 & 0 \end{pmatrix},\,\,
Z=\begin{pmatrix} i & 0 \\ 0 & -i \end{pmatrix}. $$
The invariant vector fields $X_A,X_B$ generated by $A,B\in{\frak su}(2)$ 
can be expressed in terms of the fields
$\frac{\partial}{\partial\theta}, \frac{\partial}{\partial\phi}$. 
Given $\mu\in{\cal O}_{\eta}$,
$X_B(\mu)$ is defined by the curve $e^{tB}\mu$. If $\mu=g\eta$, with $g$ as above, then
$e^{tB}g$ is the element of $SU(2)$ determined by the pair
 $$(x',\,y')=(x\cos t-\bar y\sin t,\, y\cos t +\bar x\sin t).$$  
 An easy but tedious calculation
 shows  that 
$$(x',\,y')=\big(\cos(\theta'/2)e^{i\phi_1'},\,\sin(\theta'/2)e^{-i\phi_2'}  \big)+O(t^2),$$
with 
$\theta'=\theta+2t\cos\phi,\,\,\,\phi_1'=\phi_1+t\tan(\theta/2)\sin\phi,\,\,\,
\phi_2'=\phi_2+t \cot (\theta/2)\sin\phi.$
Therefore
\begin{equation}\label{XB}
X_B(\theta,\phi)=2\cos\phi\frac{\partial}{\partial\theta}
-2\cot\theta\sin\phi\frac{\partial}{\partial\phi}
\end{equation}
Similarly
\begin{equation}\label{XA}
X_A(\theta,\phi)=2\sin\phi\frac{\partial}{\partial\theta}
+2\cot\theta\cos\phi\frac{\partial}{\partial\phi}
\end{equation}

The symplectic structure on ${\cal O}_{\eta}$ is defined by the
form $\omega$, whose action on invariant vector fields is 
$$\omega_{\mu}(X_C(\mu),X_D(\mu))=\mu([C,D]).$$
$\omega$ 
can also be expressed 
in the spherical coordinates. With the above notations
$$\omega_{\mu}(X_A,X_B)=\eta(g^{-1}[A,B]g)=-2k(|x|^2-|y|^2)=-2k\cos\theta.$$
Using  (\ref{XB})   and (\ref{XA})  a simple calculation gives
\begin{equation}\label{omegaSU2}
\omega=\frac{k}{2}\sin\theta d\theta\wedge d\phi.
\end{equation}
Given $C\in{\frak su}(2)$, the function $h_C$ on ${\cal O}_{\eta}$
defined by $h_C(\mu)=\mu(C)$
satisfies $\omega(X_C,\,.)=dh_C$.
In spherical coordinates
  $$h_A(\mu)=g\eta(A)=\eta(g^{-1}Ag)=-k(xy+\bar x\bar y)=-k\sin\theta\cos\phi.$$
That is, 
\begin{equation}\label{hASU2}
h_A(\theta,\phi)=-k\sin\theta\cos\phi.
\end{equation}
Moreover $h_A$ satisfies the normalisation condition
$\int_{S^2}h_A\omega=0.$
A similar calculation gives
\begin{equation}\label{hBSU2}
h_B(\theta,\phi)=k\sin\theta\sin\phi
\end{equation}

Henceforth we asume that $k=\frac{n}{2\pi}$, with $n\in{\Bbb Z}$. Then 
the orbit ${\cal O}_{\eta}$ possesses an invariant prequantization (see
\cite{bK70})

We can consider the family $\{\psi_t\}$ of symplectomorphisms of ${\cal O}_{\eta}$
defined by $\psi_t(\mu):=e^{tA}\cdot\mu$.  As 
\begin{equation}\label{esA}
e^{tA}=\begin{pmatrix} \cos t & i\sin t \\
       i\sin t & \cos t \end{pmatrix},
\end{equation}
 hence $\psi_{\pi}:S^2\rightarrow S^2$ is the identity, and $\psi=\{\psi_t\,|\, t\in[0,\pi]\}$
is a loop in the group of $\text{Ham}({\cal O}_{\eta})$.

If one takes the north pole $p(\theta=0,\phi=0)$, the curve $\psi_t(p)$ is the path  obtained 
as product of the paths defined by
 the meridians $\phi=\pi/2$ and $\phi=3\pi/2$. So by (\ref{hASU2})  $h_A(\psi_t(p))=0$, and 
$$S=\{(\theta,\phi)\,|\, \pi/2\leq \phi\leq 3\pi/2,\,\,\,\theta\in[0,\pi]\}$$
oriented with $d\theta\wedge d\phi$ is an oriented surface whose boundary is the curve
$\psi_t(p)$. By (\ref{omegaSU2}) $\int_S\omega =k\pi$, and  
from (\ref{kappa})
 we obtain   
$\kappa_p(\psi)=(-1)^n$.

We could calculate $\kappa_q(\psi)$ for
$q(\theta=\pi/2,\phi=0)$.  Now $\psi_t(q)=q$ for all $t$, 
hence the integral 
of $\omega$ in (\ref{kappa}) vanishes. $h_A(q)=-n/2\pi$,  consequently
$-\int_0^{\pi}f_t(\psi_t(q))=-n/2$, and $\kappa_q(\psi)=(-1)^n$.

Let us consider the point $r=(\pi/2,\pi/2)\in S^2$, according to
(\ref{xyinS2}) this point can be represented
by the element of $g\in SU(2)$ defined by $x=2^{-1/2}i,\,\, y=2^{-1/2}$. Denoting
by $(\theta',\phi')$ the spherical
coordinates of $\psi_t(r)$, from (\ref{esA}) one deduces
$$e^{i\phi'}\cos(\theta'/2)=\frac{i}{\sqrt 2}(\cos t-\sin t),\,\,\,\,\,\,
\sin(\theta'/2)=\frac{1}{\sqrt 2}(\cos t+\sin t)$$
Hence $\theta'=2t+\pi/2,\,\,\phi'=\pi/2$ when $t\in[0,\pi/4]$, etc. That is,
$\{\psi_t(r)\}$ is the union of the meridians $\phi=\pi/2$ and $\phi=3\pi/2$.
So $h_A(\psi_t(r))=0$. On the other hand
$$\int_0^{\pi} \int_{\pi/2}^{3\pi/2}\frac {n}{4\pi}\sin\theta d\theta\wedge d\phi
=\frac{n}{2}.$$
So $\kappa_r(\psi)=(-1)^n.$ 

The equalities $\kappa_p(\psi)=\kappa_q(\psi)=  \kappa_r(\psi)$     can  also be considered as
a  checking  of Theorem \ref{Therefi}.  

\smallskip

We will determine $\kappa(\chi)$, when $\chi_t$ is the symplectomorphism of $S^2$
given by  $\chi_t(q)=e^{t(aA+bB)}q$, where 
 $a,b\in{\Bbb R}$.
For $t\geq 0$ we put $c=t(b+ai)$, so
$$t(aA+bB)=\begin{pmatrix} 0 & c \\
-\bar c & 0
\end{pmatrix}.$$
The matrix ${t(aA+bB)}$ can be diagonalized, and 
 $$D^{-1}t(aA+bB)D=\text{diag}(i|c|,\,-i|c|),$$
with
 $$ D=\frac{1}{\sqrt 2}\begin{pmatrix} -i\epsilon & -1 \\
                                   1 & i\bar\epsilon \end{pmatrix},$$
where $\epsilon=c/|c|$.
Hence
$$e^{t(aA+bB)}=D\,\text{diag}(e^{i|c|},\,e^{-i|c|})\,D^{-1}.$$
It is straightforward to deduce
\begin{equation}\label{exponencial}
e^{t(aA+bB)}=\begin{pmatrix} \cos |c| & \epsilon\sin|c|  \\
-\bar\epsilon\sin|c| & \cos |c|
\end{pmatrix}.
\end{equation}
 For
$t_1=\pi/\sqrt{a^2+b^2}$ the Hamiltonian symplectomorphism $\chi_{t_1}=\text{id}$,
so
$\{\chi_t\,|\, t\in[0,t_1]  \}$
 is a loop in $\text{Ham}(S^2)$. From now on we assume $\sqrt{a^2+b^2}=1$, then
$\chi_{\pi}=\text{id}$.

Let $p$ be the north pole, then $\chi_t(p)$ is the point which
corresponds to the pair 
\begin{equation}\label{pairxy}
(x=\cos t,\, y=\epsilon\sin t)
\end{equation}
 in the notation (\ref{xysu2}).
We put $\epsilon=e^{i\alpha}$; from (\ref{xyinS2}) and  (\ref{pairxy}) it follows
that the spherical coordinates of $\chi_t(p)$ 
 are $(2t,\alpha)$, for  $t\in[0,\,\pi/2]$.

Similarly, when $t$ runs on $[\pi/2,\, \pi]$ the point
$\chi_t(p)$ runs on the meridian $\phi=\pi+\alpha$ from $\theta=\pi$ to $\theta=0$; that is,
$\chi_t(p)=(2\pi-2t,\pi+\alpha)$.

Since $h_{aA+bB}=ah_A+bh_B$, by (\ref{hASU2}) and (\ref{hBSU2})
$$h_{aA+bB}(\theta,\phi)=k\sin\theta(-a\cos\phi+b\sin\phi).$$
As $\epsilon= (a^2+b^2)^{-1/2}(b+ai)=\cos\alpha+i\sin\alpha$, we obtain
 $$h_{aA+bB}(\theta,\phi)=k\sin\theta\sin(\phi-\alpha).$$
Taking into account the spherical coordinates of $\chi_t(p)$ determined above, one deduces
 $h_{aA+bB}(\chi_t(p))=0$, for every $t\in[0,\pi]$.  
Thus $\kappa(\chi)=\exp(2\pi i\int_S\omega)$, where $S$ is the hemisphere
limited by the meridian $\phi=\alpha$ and $\phi=\pi+\alpha$. Therefore $\kappa(\chi)=(-1)^n$.
In summary
\begin{Thm}\label{kappasu2} Let $\eta$ be the element of ${\frak su}(2)^*$ defined by
$\eta\begin{pmatrix} ai & w \\ -\bar w & -ai 
\end{pmatrix}=
 \frac{n}{2\pi}a$, with $n\in{\Bbb Z}$. If $\chi$ is a loop in $\text{Ham}({\cal O}_{\eta})$
which is a $1$-parameter subgroup generated by an invariant vector field, then
 $\kappa(\chi)=(-1)^n.$
\end{Thm}

The vector $aA+bB\in{\frak su}(2)$, with $a^2+b^2=1$ can be deformed by means of a rotation
into $a'A+b'B$, if $(a')^{2}+(b')^2=1$. If we denote $\chi'_t:=\text{exp}(t(a'A+b'B))$, 
by Theorem \ref{invakappahomot}
$\kappa(\chi)=\kappa(\chi')$. Therefore Theorem \ref{kappasu2} can also be considered as a checking of 
Theorem \ref{invakappahomot}.

\end{document}